\documentclass[11pt]{amsart}

\usepackage[T1]{fontenc}
\usepackage[utf8]{inputenc}
\usepackage[margin=1.3in]{geometry}
\usepackage{amssymb}
\usepackage{microtype}
\usepackage{nicefrac}
\usepackage{esint}
\usepackage[colorlinks=true,linkcolor=blue,citecolor=blue,urlcolor=blue]{hyperref}
\usepackage{xcolor}
\usepackage{tikz}
\usetikzlibrary{arrows.meta, positioning, calc}
\makeatletter
\renewcommand{\paragraph}{\@startsection{paragraph}{4}%
  {\z@}{1.0ex \@plus .5ex \@minus .2ex}{-1em}%
  {\normalfont\normalsize\bfseries}}
\makeatother
\newcommand{\R}{\mathbb{R}}
\DeclareMathOperator*{\esssup}{ess\,sup}
\DeclareMathOperator*{\essinf}{ess\,inf}
\newcommand{\Lean}{Lean~4}
\newcommand{\mathlib}{mathlib}
\newcommand{\leanversion}{v4.29.0-rc6}
\newcommand{\mathlibcommit}{5c8398d}
\newcommand{\projectcommit}{6614c52}

\newtheorem{theorem}{Theorem}[section]
\theoremstyle{remark}

\title[Formalization of De Giorgi--Nash--Moser Theory]{Formalization of De Giorgi--Nash--Moser Theory in Lean}

\author{Scott Armstrong}
\address{CNRS \& Laboratoire Jacques-Louis Lions, Sorbonne Universit\'e and Courant Institute, School of Mathematics, Computing, and Data Science, New York University.}
\email{scottnarmstrong@gmail.com}
\thanks{S.A. acknowledges support from the European Research Council (ERC) under the European Union's Horizon Europe research and innovation programme, grant agreement number 101200828.}

\author{Julia Kempe}
\address{Courant Institute, School of Mathematics, Computing, and Data Science, New York University.}
\email{kempe@nyu.edu}
\thanks{J.K. thanks the Simons Foundation for support through the Collaborative Grant ``The Physics of Learning and Neural Computation''.}


\date{April 7, 2026}

\begin{document}

\begin{abstract}
We present a formalization in \Lean{} of the core interior De Giorgi--Nash--Moser theory for uniformly elliptic divergence-form equations with bounded measurable coefficients. The formalized results include local boundedness of weak subsolutions, the weak Harnack inequality for positive weak supersolutions, Moser's Harnack inequality for positive weak solutions, and interior H\"older regularity. This is, to our knowledge, the first machine-checked formalization of a major theorem in modern PDE theory. The development also required substantial new infrastructure for Sobolev spaces on bounded domains, weak solutions of elliptic equations, and quantitative regularity estimates. More broadly, it suggests that large-scale autoformalization of hard analysis in \Lean{} is now within reach.
\end{abstract}

\maketitle

\section{Introduction}

The De Giorgi--Nash continuity estimate, proved independently by De~Giorgi~\cite{DeGiorgi} and Nash~\cite{Nash}, states that every weak solution $u \in H^1(B_1)$ of an elliptic equation of the form
\begin{equation}
-\nabla \cdot \mathbf{a} \nabla u = 0 \quad \mbox{in} \ B_1 \subseteq\mathbb{R}^d\,,
\label{e.pde}
\end{equation}
where~$B_1$ is the unit ball in~$\mathbb{R}^d$ and~$\mathbf{a}$ is a measurable function from~$B_1$ to the set of invertible, real~$d\times d$ matrices satisfying, for some~$0< \lambda \leq \Lambda<\infty$, the uniform ellipticity condition
\begin{equation}
\xi \cdot \mathbf{a}(x) \xi \geq \lambda |\xi|^2
\quad \mbox{and} \quad
\xi \cdot \mathbf{a}^{-1}(x) \xi \geq \Lambda^{-1} |\xi|^2
\quad \mbox{for a.e.} \ x \in B_1, \ \xi \in \mathbb{R}^d\,,
\label{e.ue}
\end{equation}
is $\alpha$-H\"older continuous in~$B_{\nicefrac12}$ (more precisely, there exists a version of~$u$ which is H\"older continuous) for an exponent~$\alpha(d,\Lambda/\lambda) \in (0,\nicefrac12]$, and satisfies the estimate
\begin{equation*}
\| u \|_{L^\infty(B_{\nicefrac12})}
+
\sup_{x,y\in B_{\nicefrac12}, \, x\neq y}
\frac{|u(x) - u(y)|}{|x-y|^\alpha}
\leq
C\| u \|_{L^2(B_{1})}\,.
\end{equation*}
The prior state of the art, rooted in Schauder theory, required H\"older regularity of the coefficients to obtain regularity of solutions. It was a remarkable discovery that the purely measure-theoretic coercivity condition~\eqref{e.ue} is sufficient to force continuity of solutions, a pointwise property. This conceptual and technical breakthrough resolved Hilbert's nineteenth problem on the regularity of minimizers of uniformly convex variational problems. Several years later, J\"urgen Moser~\cite{Moser60, Moser61} gave a third proof of H\"older regularity and established the stronger statement that nonnegative solutions of~\eqref{e.pde} satisfy a Harnack inequality, namely 
\begin{equation}
\esssup_{B_{\nicefrac12}} u \leq C \essinf_{B_{\nicefrac12}} u\,.
\label{e.Harnack}
\end{equation}
These estimates and the circle of ideas surrounding them are known as \emph{De Giorgi--Nash--Moser theory} and are now foundational across elliptic and parabolic PDE, geometric analysis, and probability theory.

\smallskip

This paper describes a complete formalization using the \Lean{} proof assistant~\cite{Lean4} of the core interior De Giorgi--Nash--Moser theory for uniformly elliptic divergence-form equations. We formalize local boundedness of weak subsolutions, the weak Harnack inequality for positive weak supersolutions, Moser's Harnack inequality for positive weak solutions, and interior H\"older regularity of solutions. The formalized theorem statements are given in Section~\ref{sec:statements}, and the classical proof route they follow is outlined informally in Section~\ref{sec:outline}.

This project required developing a collection of analytic tools that were not yet available in \Lean, and in some cases had not previously been formalized in the form needed here. The proofs of De Giorgi--Nash--Moser theory intertwine weak derivatives, Sobolev spaces on bounded domains, functional inequalities for Sobolev functions, truncation arguments, and nonlinear test functions. A substantial part of the work therefore lay in building a reusable infrastructure for weak solutions of elliptic equations before the regularity theory itself could be carried out.

Against this background, the present work enters genuinely new territory. To our knowledge, it is the first formalization in \Lean{} of a nontrivial theorem about solutions of a partial differential equation, and in particular the first Lean formalization of a regularity result for weak solutions of a PDE. The closest directly relevant Lean antecedents are the recent formalization of the Gagliardo--Nirenberg--Sobolev inequality by van Doorn and Macbeth~\cite{vanDoornMacbeth}, and Doll's formalization of Schwartz functions, tempered distributions, and Sobolev spaces on $\R^d$ defined via the Fourier transform~\cite{Doll}. The $h$-principle formalized by Massot, van Doorn, and Nash~\cite{MassotVanDoornNash} concerns partial differential relations rather than the analytic theory of elliptic equations. Beyond Lean, the Coq formalization of the Lax--Milgram theorem by Boldo, Cl\'ement, Faissole, Martin, and Mayero~\cite{LaxMilgramCoq} provides an abstract existence and uniqueness tool for weak problems. These earlier achievements are important precursors, but they address either foundational infrastructure or substantially simpler PDE settings, rather than a full regularity theory for weak solutions of rough-coefficient elliptic equations. At the same time, our development relies crucially on the broader analytic library already present in \mathlib, especially the formalized Gagliardo--Nirenberg--Sobolev inequality and the surrounding whole-space measure-theoretic and functional-analytic tools.

\smallskip

The Lean code is available at
\begin{center}
\url{https://github.com/scottnarmstrong/DeGiorgi}
\end{center}
The version described in this paper corresponds to commit \texttt{\projectcommit}; it is built with \Lean{}~\leanversion{} and \mathlib{} commit \texttt{\mathlibcommit}. The repository also contains the file \texttt{DeGiorgi\_Lean\_to\_Tex.tex}, a machine-generated natural-language translation of the Lean development.

\smallskip

The formal development was carried out with extensive use of large language models (LLMs), under close supervision by the authors. We supplied detailed proof blueprints, intervened whenever the formal development stalled or required a reorganization of the argument, and validated the final mathematical claims against the checked Lean artifact.

\smallskip

In the concluding section of~\cite{Doll}, Doll wrote that ``it is currently not possible to prove [formalize] non-trivial results about solutions of partial differential equations.'' Our experience in this project is that recent advances in large language models (LLMs) have changed the situation substantially, making a formalization of this scale practically feasible in Lean.

The paper is organized as follows. Section~\ref{sec:statements} states the formalized theorem package in ordinary mathematical language. Section~\ref{sec:outline} then sketches the informal proof architecture. Section~\ref{sec:infrastructure} describes the analytic infrastructure that had to be built in order to support the proof. Section~\ref{sec:challenges} discusses the main formalization challenges, and Section~\ref{sec:conclusion} concludes.

\section{Formalized Statements of the Main Results}
\label{sec:statements}

Throughout, we work on the ambient Euclidean space~$\R^d$ for a fixed dimension~$d \geq 3$, and we write~$B_r$ for the open ball~$B(0,r)$ centered at the origin. This restriction entails no real loss of generality: in dimension~$d = 2$, a simpler argument already gives H\"older regularity, and the two-dimensional case may also be deduced from the three-dimensional one by adding a dummy variable. We work in dimension $d \ge 3$ because the Sobolev embedding $W^{1,2} \hookrightarrow L^{2^*}$ with $2^* = 2d/(d-2)$ is then supercritical and including $d=2$ would require splitting into a separate case.  

\subsection{Analytic setting}

\paragraph{Sobolev spaces.}
\label{sec:sobolev}
Let $\Omega \subseteq \R^d$ be open. A measurable function $g \colon \Omega \to \R$ is a weak partial derivative of $f \colon \Omega \to \R$ in the $i$-th coordinate direction if
\[
  \int_\Omega f(x)\,\partial_i \varphi(x)\,dx
  =
  -\int_\Omega g(x)\,\varphi(x)\,dx
\]
for every test function $\varphi \in C_c^\infty(\Omega)$. We write $f \in W^{1,p}(\Omega)$ if $f \in L^p(\Omega)$ and each weak partial derivative belongs to $L^p(\Omega)$. The space $W^{1,2}_0(\Omega)$ is defined as the closure of $C_c^\infty(\Omega)$ in the $W^{1,2}$ topology.

\paragraph{Elliptic coefficients.}
\label{sec:coefficients}
An elliptic coefficient field on $\Omega$ is a measurable map $\mathbf{a} \colon \Omega \to \R^{d \times d}$ together with positive constants $0 < \lambda \le \Lambda < \infty$ such that, for Lebesgue-a.e.\ $x \in \Omega$ and every $\xi \in \R^d$,
\begin{equation}
\xi \cdot \mathbf{a}(x)\xi \geq \lambda |\xi|^2
\qquad \text{and} \qquad
\xi \cdot \mathbf{a}(x)^{-1}\xi \geq \Lambda^{-1} |\xi|^2.
\label{e.coercive}
\end{equation}
For most of the final theorem statements it is convenient to normalize to $\lambda = 1$, so that the dependence of the estimates is displayed only through $\Lambda$.

\paragraph{Weak solutions.}
\label{sec:weak-solutions}
Given an elliptic coefficient field $\mathbf{a}$ on $\Omega$ and functions $u,\varphi \in W^{1,2}(\Omega)$, we define
\[
  B_{\mathbf{a}}(u,\varphi)
  :=
  \int_\Omega \nabla \varphi(x) \cdot \mathbf{a}(x)\nabla u(x)\,dx.
\]
We say that $u$ is a weak subsolution of $-\nabla \cdot (\mathbf{a}\nabla u) \le 0$ on $\Omega$ if $u \in W^{1,2}(\Omega)$ and
\[
  B_{\mathbf{a}}(u,\varphi) \le 0
\]
for every nonnegative $\varphi \in W^{1,2}_0(\Omega)$. Weak supersolutions and weak solutions are defined analogously.

\subsection{Main theorems}

\begin{theorem}[Local boundedness of subsolutions]\label{t.linfty}
Let $d \geq 3$ and let $\mathbf{a}$ be an elliptic coefficient field on $B_1$. If $u$ is a weak subsolution of $-\nabla \cdot (\mathbf{a}\nabla u) \leq 0$ on $B_1$ and $u_+ = \max(u,0)$ belongs to $L^2(B_1)$, then for almost every $x \in B_{1/2}$,
\begin{equation}
  \max(u(x),0)
  \le
  C_{\mathrm{DG}}(d,\mathbf{a})
  \Bigl(\int_{B_1} |u_+|^2\Bigr)^{\nicefrac12}.
\label{e.linfty}
\end{equation}
\end{theorem}

\begin{theorem}[Weak Harnack inequality]\label{t.weak-harnack}
Let $d \geq 3$, let $\mathbf{a}$ be a normalized elliptic coefficient field on $B_1$, let $0 < q < 1$, and let $u$ be a positive weak supersolution of $-\nabla \cdot (\mathbf{a}\nabla u) \ge 0$ on $B_1$. Writing
\[
q^* := \frac{qd}{d-2},
\]
one has
\[
\Bigl(\int_{B_{1/4}} |u|^{q^*}\Bigr)^{1/q^*}
\le
\Bigl(\frac{C_{\mathrm{wH}}(d)}{(1-q)^{\gamma_d}}\Bigr)^{\Lambda^{\nicefrac12}}
\essinf_{B_{1/4}} u,
\]
where $C_{\mathrm{wH}}(d)$ and $\gamma_d$ depend only on $d$.
\end{theorem}

\begin{theorem}[Crossover estimate]\label{t.crossover}
Let $d \geq 3$, let $\mathbf{a}$ be a normalized elliptic coefficient field on $B_1$, and let $u$ be a positive weak supersolution on $B_1$. Then
\[
\Bigl(\fint_{B_{\nicefrac12}} |u|^c\Bigr)
\Bigl(\fint_{B_{\nicefrac12}} |u|^{-c}\Bigr)
\leq C'_{\mathrm{cross}}(d),
\qquad
c := \frac{c'_{\mathrm{cross}}(d)}{\Lambda^{\nicefrac12}},
\]
for dimension-dependent constants $c'_{\mathrm{cross}}(d)$ and $C'_{\mathrm{cross}}(d)$.
\end{theorem}

\begin{theorem}[Harnack inequality]\label{t.harnack}
Let $d \geq 3$ and let $\mathbf{a}$ be a normalized elliptic coefficient field on $B_1$. If $u$ is a positive weak solution of $-\nabla \cdot (\mathbf{a}\nabla u) = 0$ on $B_1$, then
\begin{equation}
  \esssup_{B_{1/2}} u
  \le
  \exp\bigl(C_{\mathrm{H}}(d)\Lambda^{\nicefrac12}\bigr)
  \essinf_{B_{1/2}} u.
\label{e.harnack}
\end{equation}
\end{theorem}

\begin{theorem}[Interior H\"older regularity]\label{t.holder}
Let $d \geq 3$, let $\mathbf{a}$ be a normalized elliptic coefficient field on $B_1$, let $p_0 > 1$, and let $u$ be a weak solution of $-\nabla \cdot (\mathbf{a}\nabla u) = 0$ on $B_1$ with $|u|^{p_0} \in L^1(B_1)$. Then there exists a function $v \colon \R^d \to \R$ satisfying $v = u$ almost everywhere on $B_1$, and there exists an exponent $\alpha > 0$ with
\begin{equation}
  \alpha \ge \exp\bigl(-C_{\mathrm{H\ddot{o}l}}(d)\Lambda^{\nicefrac12}\bigr),
\label{e.holder-exponent}
\end{equation}
such that, for all $x,y \in B_{1/2}$,
\begin{equation}
  |v(x)-v(y)|
  \le
  C_{\mathrm{H\ddot{o}l}}(d)\Lambda^{d/(2p_0)}
  \Bigl(\frac{p_0}{p_0-1}\Bigr)^{d/p_0}
  \Bigl(\int_{B_1} |u|^{p_0}\Bigr)^{1/p_0}
  |x-y|^\alpha.
\label{e.holder}
\end{equation}
\end{theorem}

\subsection{Quantitative form of the results}

Every constant appearing in Theorems~\ref{t.linfty}--\ref{t.holder} is a closed-form arithmetic expression in the dimension, with the dependence on ellipticity displayed explicitly in the statement. The proofs contain no compactness or contradiction arguments that hide constants behind qualitative existence statements: the bounds are assembled from named quantities by elementary operations.

For the local boundedness theorem, the dependence on the ellipticity ratio is polynomial and comes from the one-step energy estimate and the De Giorgi recurrence threshold. For Harnack, the formalization follows Bombieri--Giusti~\cite{BombieriGiusti} and recovers the sharp dependence
\[
\exp\bigl(C_{\mathrm{H}}(d)\Lambda^{\nicefrac12}\bigr)
\]
on the ellipticity parameter. The same crossover mechanism yields the lower bound~\eqref{e.holder-exponent} for the H\"older exponent. This lower bound is not expected to be sharp in $\Lambda$; however, in dimensions $d>2$ no better quantitative lower bound is currently known in the literature.

\section{Outline of the Informal Proof}
\label{sec:outline}

Having stated the formalized theorem package, we now sketch the mathematical proof route formalized in \Lean. The proof proceeds through four layers: a bounded-domain Sobolev library, nonlinear iteration schemes for sub- and supersolutions, the weak and strong Harnack inequalities, and finally oscillation decay leading to H\"older regularity.

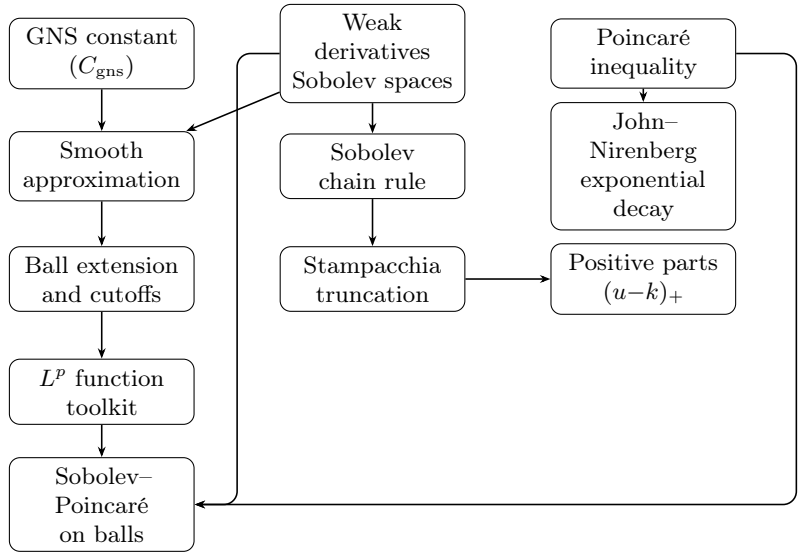
\begin{figure}[ht]
\centering
\begin{tikzpicture}[
    box/.style={draw, rounded corners, font=\footnotesize,
      text width=2.2cm, align=center, minimum height=0.55cm},
    arr/.style={-{Stealth[length=4pt]}, semithick},
  ]
  \node[box] (GNS)     at (0,0)      {GNS constant\\($C_{\mathrm{gns}}$)};
  \node[box] (Sob)     at (3.6,0)    {Weak derivatives\\Sobolev spaces};
  \node[box] (poinc)   at (7.2,0)    {Poincar\'e\\inequality};

  \node[box] (approx)  at (0,-1.5)   {Smooth\\approximation};
  \node[box] (chain)   at (3.6,-1.5) {Sobolev\\chain rule};
  \node[box] (JN)      at (7.2,-1.5) {John--Nirenberg\\exponential decay};

  \node[box] (ext)     at (0,-3.0)   {Ball extension\\and cutoffs};
  \node[box] (stamp)   at (3.6,-3.0) {Stampacchia\\truncation};
  \node[box] (pospart) at (7.2,-3.0) {Positive parts\\$(u{-}k)_+$};

  \node[box] (Lp)      at (0,-4.5)   {$L^p$ function\\toolkit};

  \node[box] (SP)      at (0,-6.0)   {Sobolev--Poincar\'e\\on balls};

  \draw[arr] (GNS) -- (approx);
  \draw[arr] (Sob) -- (approx);
  \draw[arr] (Sob) -- (chain);
  \draw[arr] (poinc) -- (JN);
  \draw[arr] (approx) -- (ext);
  \draw[arr] (chain) -- (stamp);
  \draw[arr] (stamp) -- (pospart);
  \draw[arr] (ext) -- (Lp);
  \draw[arr] (Lp) -- (SP);
  \draw[arr, rounded corners=5pt]
    (poinc.east) -- ++(0.8,0) -- ++(0,-6.0) -- (SP.east);
  \draw[arr, rounded corners=5pt]
    (Sob.west) -- (1.8,0) -- (1.8,-6.0) -- (SP.east);
\end{tikzpicture}
\caption{Sobolev infrastructure. All boxes are reusable library components not specific to De Giorgi--Nash--Moser theory.}
\label{fig:sobolev}
\end{figure}

\begin{figure}[ht]
\centering
\begin{tikzpicture}[
    box/.style={draw, rounded corners, font=\footnotesize,
      text width=2.5cm, align=center, minimum height=0.55cm},
    gbox/.style={box, fill=gray!15},
    thm/.style={box, line width=0.8pt},
    arr/.style={-{Stealth[length=4pt]}, semithick},
  ]
  \node[gbox] (sobolev) at (0, 0)   {Sobolev infrastructure\\(Fig.~\ref{fig:sobolev})};
  \node[box]  (coeff)   at (3, 0)   {Elliptic coefficients and \\weak formulation};
  \node[box]  (loggrad) at (6, 0)   {Log-gradient\\bound};
  \node[gbox] (JN)      at (9, 0)   {John--Nirenberg\\(Fig.~\ref{fig:sobolev})};

  \node[box] (cacc)    at (0, -1.7)  {Caccioppoli\\inequality};
  \node[thm] (DG)      at (3, -1.7)  {De Giorgi $L^\infty$\\(Thm~\ref{t.linfty})};
  \node[box] (cross)   at (6, -1.7)  {Crossover\\estimate};

  \node[box] (Moser)   at (3, -3.4)  {Moser iteration\\$L^p \to L^\infty$};
  \node[box] (super)   at (6, -3.4)  {Supersolution\\estimates};

  \node[box] (wH)      at (6, -5.1)  {Weak Harnack\\inequality};

  \node[thm] (harnack) at (4.5, -6.6) {Harnack inequality\\(Thm~\ref{t.harnack})};
  \node[thm] (holder)  at (4.5, -8.0) {H\"older regularity\\(Thm~\ref{t.holder})};

  \draw[arr] (sobolev) -- (cacc);
  \draw[arr] (sobolev) -- (DG);
  \draw[arr] (coeff) -- (DG);
  \draw[arr] (coeff) -- (cacc);
  \draw[arr] (coeff) -- (loggrad);
  \draw[arr] (loggrad) -- (cross);
  \draw[arr] (JN) -- (cross);
  \draw[arr] (cacc) -- (Moser);
  \draw[arr] (DG) -- (Moser);
  \draw[arr] (cross) -- (super);
  \draw[arr] (Moser) -- (super);
  \draw[arr] (super) -- (wH);
  \draw[arr] (Moser) -- (harnack);
  \draw[arr] (wH) -- (harnack);
  \draw[arr] (harnack) -- (holder);
\end{tikzpicture}
\caption{PDE proof chain. Gray boxes are Sobolev infrastructure from Figure~\ref{fig:sobolev}. The De Giorgi $L^\infty$ bound (Theorem~\ref{t.linfty}) anchors the Moser iteration at $p=2$. The crossover estimate uses the log-gradient bound and John--Nirenberg to bridge positive and negative moments. Harnack (Theorem~\ref{t.harnack}) combines the Moser bound with weak Harnack via a chain argument. H\"older regularity (Theorem~\ref{t.holder}) follows by oscillation decay.}
\label{fig:pde}
\end{figure}
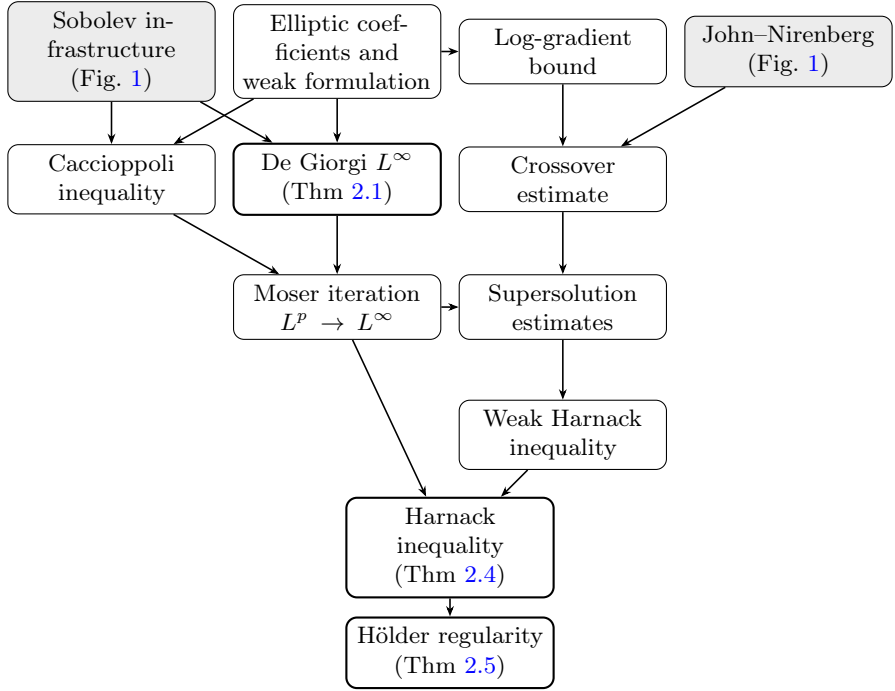

\subsection{Analytic infrastructure}

The proof begins with a bounded-domain Sobolev library. Weak derivatives are defined by integration by parts against smooth compactly supported test functions, and Sobolev spaces on an open set $\Omega \subseteq \R^d$ are defined through the existence of weak partial derivatives in $L^p(\Omega)$. Since later arguments repeatedly reuse a specific weak gradient, the formal development is organized around explicit witnesses for weak derivatives rather than bare existential statements.

Smooth approximation on balls is the gateway to the Sobolev chain rule. Once one knows that a Sobolev function can be approximated in both function and gradient by smooth compactly supported functions, one may apply the classical chain rule to the smooth approximants and pass to the limit. Combined with Stampacchia's truncation lemma, this yields the positive-part closure $u \mapsto (u-k)_+$ and makes the truncated test functions used in De Giorgi and Moser iteration admissible.

To import whole-space analytic tools from \mathlib, the formalization extends a function on $B_1$ to a compactly supported function on $\R^d$ by an explicit inversion-based extension operator. This feeds into Poincar\'e and Sobolev--Poincar\'e inequalities on balls. Finally, the John--Nirenberg theorem supplies exponential integrability from BMO bounds, which is the key bridge between logarithmic Caccioppoli estimates and the crossover argument.

\subsection{De Giorgi and Moser iterations}
\label{sec:iterations}

The $L^\infty$ bound for subsolutions is obtained by two related iteration schemes. The \emph{De Giorgi iteration} works with truncation levels. Given a weak subsolution~$u$ and a level $k > 0$, one inserts the test function $\eta^2(u-k)_+$, where $\eta$ is a cutoff between concentric balls, into the subsolution inequality to obtain a Caccioppoli estimate. Combined with Sobolev--Poincar\'e, this yields a nonlinear recurrence
\[
  Y_{n+1} \leq C_{\mathrm{rec}} B^n Y_n^{1+\alpha}
\]
for the energies $Y_n = \int_{B_{r_n}} |(u-k_n)_+|^2$ on shrinking balls and rising levels. A general recurrence closeout lemma then shows that $Y_n \to 0$ once the initial energy is below an explicit threshold, giving local boundedness.

The \emph{Moser iteration} works instead with powers. Testing with $\eta^2 |u_+|^{p-1}$ promotes $L^p$ control to $L^{p^*}$ control, where $p^* = pd/(d-2)$, and iteration along a geometric sequence of exponents gives an $L^p \to L^\infty$ bound. In our formal development, the De Giorgi $L^2 \to L^\infty$ estimate provides the qualitative boundedness needed to anchor the Moser argument at $p=2$; it is the Moser version, with its flexibility in the starting exponent, that feeds the Harnack and H\"older theorems downstream.

We decided to formalize both iteration arguments. In our development they are not logically independent, since the upper bound provided by De Giorgi's iteration is used as a qualitative input in the Moser iteration. This allows us to avoid the usual truncation step needed to produce an admissible test function at the start of the Moser argument.

\subsection{Weak Harnack and the crossover estimate}

The weak Harnack inequality controls an $L^q$ norm of a positive supersolution from below by its essential infimum. Its proof combines three ingredients.

The first is a forward estimate for positive supersolutions, obtained by applying Moser iteration to negative powers of the solution. The second is a logarithmic Caccioppoli estimate: after a suitable regularization of the logarithm, inserting $\varphi^2/u$ into the supersolution inequality gives
\[
  \int \varphi^2 \frac{|\nabla u|^2}{u^2} \leq \frac{4\Lambda}{\lambda} \int |\nabla\varphi|^2.
\]
Writing $v = -\log u$, this gives control of $\nabla v$ in $L^2$ and hence, via Poincar\'e, control of the BMO (bounded mean oscillation) norm of $v$, namely
\[
  \|v\|_{\mathrm{BMO}(B)}
  := \sup_{B' \subseteq B} \fint_{B'} |v - (v)_{B'}|,
  \qquad
  (v)_{B'} := \fint_{B'} v.
\]

The third ingredient is the Bombieri--Giusti crossover estimate. John--Nirenberg upgrades the BMO control of $v$ to exponential integrability, which in turn bounds the product of positive and negative moments of $u$:
\[
  \Bigl(\fint_{B_{\nicefrac12}} |u|^c \Bigr) \Bigl(\fint_{B_{\nicefrac12}} |u|^{-c}\Bigr)
  \leq C'(d),
  \qquad
  c = c'(d)\Lambda^{-\nicefrac12}.
\]
This exponent is what ultimately produces the sharp $\exp(C\Lambda^{\nicefrac12})$ dependence in the Harnack constant.

\subsection{From local estimates to regularity}

The Harnack inequality (Theorem~\ref{t.harnack}) is assembled from the Moser $L^p \to L^\infty$ bound and the weak Harnack inequality through a chain argument. The half-ball $B_{\nicefrac12}$ is covered by 17 overlapping quarter-balls, each centered in $B_{\nicefrac12}$, such that consecutive balls overlap by a definite amount. Multiplying the local Harnack ratios along the chain gives the global estimate
\[
  \esssup_{B_{1/2}} u \leq \exp\bigl(C_{\mathrm{H}}(d)\Lambda^{\nicefrac12}\bigr)\essinf_{B_{1/2}} u.
\]

Interior H\"older regularity (Theorem~\ref{t.holder}) then follows from Harnack by the standard oscillation-decay argument. Applying Harnack to $\esssup_{B_r}u-u$ and to $u-\essinf_{B_r}u$ shows that the oscillation decreases by a fixed factor when the radius is halved. Iterating this decay yields a H\"older modulus of continuity, and the formalization constructs a continuous representative directly from nested essential suprema and infima on dyadic balls.

\section{Formalizing the Analytic Infrastructure}
\label{sec:infrastructure}

A substantial part of the project consisted of building the analytic infrastructure needed before the De Giorgi--Nash--Moser arguments themselves could begin. The Sobolev library alone comprises roughly 20{,}000 lines of Lean and is, to our knowledge, the first formalization in any proof assistant of Sobolev spaces on bounded domains via weak derivatives. This is complementary to Doll's Fourier-analytic construction of Sobolev spaces on $\R^d$~\cite{Doll}.

\subsection{Weak derivatives and Sobolev spaces}

Weak partial derivatives are defined by integration by parts against smooth compactly supported test functions, and $W^{1,p}(\Omega)$ is defined by combining $L^p$ membership of the function with $L^p$ membership of its weak partial derivatives. In the formal development, however, existential statements about weak gradients are rarely sufficient for downstream use. Many later arguments need to reuse the same chosen gradient repeatedly, so the library is organized around explicit witness structures that package a function together with a specific weak gradient and its $L^p$ bounds.

The space $W^{1,2}_0(\Omega)$ is defined by approximation in both function and gradient by smooth compactly supported functions. This is the formulation best suited to the PDE applications, since the admissible test functions in the weak formulation naturally live in $W^{1,2}_0$. It is also the formulation that interacts most cleanly with the chain rule, truncation arguments, and the extension operator discussed below.

\subsection{Approximation, truncation, and extension}

The Sobolev chain rule is proved by approximation: if $u \in W^{1,2}(\Omega)$ and $\Phi \colon \R \to \R$ is smooth with bounded derivative and $\Phi(0)=0$, then $\Phi \circ u \in W^{1,2}(\Omega)$ with weak gradient $\Phi'(u)\nabla u$. The proof passes from smooth approximants to the limit in both function and gradient. This is one of the key interfaces between the classical differential calculus available in \mathlib{} and the weak-derivative framework needed for PDE.

Combined with Stampacchia's truncation lemma, the chain rule gives the positive-part closure
\[
(u-k)_+ \in W^{1,2}(\Omega),
\qquad
\nabla (u-k)_+ = \nabla u\,\mathbf{1}_{\{u>k\}},
\]
which is the analytic fact that legitimizes the truncated test functions used throughout the De Giorgi and Moser iterations.

To pass from balls to the whole-space inequalities available in \mathlib, the formalization constructs an explicit extension operator from $B_1$ to $\R^d$ using sphere inversion on the shell $1<|x|<2$ together with a radial cutoff. This operator preserves Sobolev regularity with quantitative control. Although the idea is classical, verifying it at the level of weak derivatives, approximation, and transport of witnesses is one of the longest pieces of the supporting library.

\subsection{Poincar\'e, Sobolev--Poincar\'e, and John--Nirenberg}

The Poincar\'e inequality on balls is proved from a representation formula obtained by integrating the fundamental theorem of calculus along rays. After subtracting the average and extending to the whole space, one obtains Sobolev--Poincar\'e on balls by applying the Gagliardo--Nirenberg--Sobolev inequality from \mathlib{} and restricting back. These inequalities are the quantitative engine behind both the De Giorgi and Moser iterations.

The John--Nirenberg theorem provides the other major analytic input. Once one has a BMO bound for $\log u$ on a ball, John--Nirenberg upgrades it to exponential integrability, and this is exactly what drives the crossover estimate. From the point of view of the formalization, this is one of the places where a substantial general-purpose analytic library pays off directly inside a delicate PDE argument.

\section{Main Formalization Challenges}
\label{sec:challenges}

The De Giorgi--Nash--Moser theory has been textbook mathematics for decades. The interesting formalization issues therefore do not arise from any mystery in the mathematics, but from the gap between the way these arguments are written on paper and the way they must be organized for a proof assistant. A recurring theme is that textbook proofs are often existential and modular only at an informal level: one says that a weak gradient exists, that a cutoff can be chosen, or that a function may be transported to another ball without naming the associated data. In the formal development, this level of abstraction is insufficient for downstream reuse.

\subsection{\texorpdfstring{Type-class synthesis and quotient $L^p$ spaces}{Type-class synthesis and quotient Lp spaces}}
\label{sec:lp-challenge}

One of the most time-consuming obstacles in the project was not mathematical but computational. In \mathlib, $L^p$ spaces are defined as quotients of measurable functions by almost-everywhere equality. Each time Lean is asked to move between a concrete function $f \colon \R^d \to \R$ and an element of this quotient space, it must resolve a substantial chain of type-class instances. For the types appearing in our development, these compatibility checks can exceed Lean's internal resource budget.

Our workaround was to avoid the abstract quotient space whenever possible. Instead, we built a thin bare-function toolkit that establishes the completeness and convergence results needed in the project directly for concrete functions equipped with the $L^p$ seminorm. This avoids the expensive elaboration path through the quotient construction while preserving exactly the analytic results needed for limit passages in the Sobolev and PDE arguments.

\subsection{Truncated test functions}
\label{sec:stampacchia-challenge}

In a textbook proof, one inserts $\eta^2 (u-k)_+$ into the weak formulation and moves on. Formalizing this step required a full infrastructure for the Sobolev chain rule, Stampacchia's truncation lemma, and positive-part closure in $W^{1,2}$. The chain rule itself already involves a delicate limit passage from smooth approximants; Stampacchia's lemma reduces to a one-dimensional statement about absolutely continuous functions and then returns to higher dimensions through Fubini. Only after these ingredients are in place can one certify that the standard truncated test functions are actually admissible in the weak formulation.

\subsection{Extension, rescaling, and change of variables}
\label{sec:extension-challenge}

Another major challenge was the extension-and-transport layer connecting local Sobolev arguments on balls to the whole-space tools available in \mathlib. In informal mathematics, one moves freely between a function on $B_1$, an extension to $\R^d$, a rescaled copy on another ball, and a smooth approximation. In Lean, each of these operations had to be implemented as a theorem with the correct hypotheses, transported weak-gradient data, and quantitative norm bounds.

The extension operator itself is explicit and classical, but proving that it preserves Sobolev regularity for rough functions required a substantial amount of geometry and analysis: inversion formulas, change-of-variables estimates, approximation by smooth functions, and the transport of chosen weak-gradient witnesses through restriction, extension, and rescaling.

\subsection{The crossover estimate}
\label{sec:crossover}

The crossover estimate is the most delicate step in the proof of the weak Harnack inequality. Mathematically, the argument we formalized is the Bombieri--Giusti one: write $v=-\log u$, obtain a BMO bound from a logarithmic Caccioppoli estimate, and then apply John--Nirenberg to control positive and negative moments simultaneously. Formally, however, the proof has to manage several layers of regularization, since neither $\log u$ nor the test function $\varphi^2/u$ is directly available in the Sobolev framework.

The solution is to replace the logarithm by smooth regularizations with bounded derivative, prove the energy estimates at the regularized level, and then pass to the limit. A second regularization is used to replace $u$ by $u+\varepsilon$ in the crossover product so that negative powers are legitimate. Much of the formal work lies in coordinating these regularizations with the Poincar\'e and John--Nirenberg machinery and in making the limit passages compatible with the chosen weak-gradient witnesses.

\subsection{Tracking explicit constants}

A design goal of the project was that the final theorem statements should display the dependence on ellipticity explicitly, rather than hide it behind qualitative constants depending on $(d,\Lambda/\lambda)$. This required naming and propagating constants throughout the entire development. In particular, the passage from the crossover exponent $c'(d)\Lambda^{-1/2}$ to the Bombieri--Giusti Harnack constant $\exp(C(d)\Lambda^{1/2})$ had to be formalized without any informal ``absorption of constants'' arguments.

The resulting quantitative bookkeeping is one of the mathematical virtues of the formalization. Every constant in the final chain is a named definition with an explicit body, and the Harnack and H\"older estimates display their $\Lambda$-dependence directly in the conclusion.

\section{Conclusion and Future Work}
\label{sec:conclusion}

We have described a complete formalization in \Lean{} of the core interior De Giorgi--Nash--Moser theory on Euclidean balls: local boundedness of weak subsolutions, the crossover estimate and weak Harnack inequality for positive supersolutions, and Harnack and interior H\"older regularity for weak solutions. Beyond the headline theorems, the project contributes a reusable analytic infrastructure for weak derivatives, Sobolev spaces on bounded domains, approximation and extension arguments, Poincar\'e and Sobolev--Poincar\'e inequalities, truncation lemmas, and John--Nirenberg theory.

Methodologically, the project was carried out with extensive use of large language models under close supervision by the authors. The present project suggests that LLM-assisted workflows can materially accelerate long formalizations in hard analysis when they are paired with explicit proof planning, careful interface design, and rigorous validation against a small trusted kernel.

\bibliographystyle{amsalpha}
\bibliography{references}

\end{document}